\documentclass[12pt]{article}
\usepackage{amssymb}
\usepackage{latexsym,bm}
\usepackage{graphicx}
\usepackage{amsmath}
\usepackage{mathrsfs}

\newcommand{\bea}{\begin{eqnarray*}}
\newcommand{\eea}{\end{eqnarray*}}
\newcommand{\be}{\begin{equation}}
\newcommand{\ee}{\end{equation}}
\newcommand{\ben}{\begin{eqnarray*}}
\newcommand{\een}{\end{eqnarray*}}

\date{}
\begin{document}
\title{Possible cardinalities of the center of a graph\footnote{E-mail addresses:
{\tt huyanan530@163.com}(Y.Hu),
{\tt zhan@math.ecnu.edu.cn}(X.Zhan).}}
\author{\hskip -10mm Yanan Hu and Xingzhi Zhan\thanks{Corresponding author.}\\
{\hskip -10mm \small Department of Mathematics, East China Normal University, Shanghai 200241, China}}\maketitle
\begin{abstract}
 A central vertex of a graph is a vertex whose eccentricity equals the radius. The center of a graph is the set of all central vertices. The central ratio
 of a graph is the ratio of the cardinality of its center to its order. In 1982, Buckley proved that every positive rational number not exceeding
 one is the central ratio of some graph. In this paper, we obtain more detailed information by determining which cardinalities are possible for the center of a graph with given order and radius. There are unexpected phenomena in the results. For example, there exists a graph of order $14$ and radius $6$ whose center has cardinality $s$ if and only if $s\in \{ 1, 2, 3, 4, 9,10,11,12,14\}.$  We also prove a related uniqueness result.
\end{abstract}



\section{Introduction}

We consider finite simple graphs. The {\it order} of a graph is its number of vertices. We denote by $V(G)$ and $E(G)$ the vertex set
 and edge set of a graph $G$ respectively. Denote by $d_{G}(u,v)$ the distance between two vertices $u$ and $v$ in $G.$ If the graph $G$ is clear from the context, we simply write $d(u,v).$
 The {\it eccentricity}, denoted by $e(v),$ of a vertex $v$ in a graph $G$ is the distance to a vertex farthest
from $v.$ Thus $e(v)={\rm max}\{d(v,u)| u\in V(G)\}.$ If $e(v)=d(v,x),$ then the vertex $x$ is called an {\it eccentric vertex} of $v.$
The {\it radius} of a graph $G,$ denoted ${\rm rad}(G),$ is the minimum eccentricity of all the vertices in $V(G),$ whereas the {\it diameter}
of $G,$ denoted ${\rm diam}(G),$ is the maximum eccentricity. A vertex $v$ is a {\it central vertex} of $G$ if $e(v)={\rm rad}(G).$ The {\it center} of a graph $G,$ denoted $C(G),$ is the set of all central vertices of $G.$ A vertex $u$ is a {\it peripheral vertex} of $G$ if $e(u)={\rm diam}(G).$ A graph with a finite radius or diameter is necessarily connected. If ${\rm rad}(G)={\rm diam}(G),$ then the graph $G$ is called {\it self-centered.} Thus, a self-centered graph
is a graph in which every vertex is a central vertex. A graph is called {\it trivial} if it is of order $1;$ otherwise it is {\it nontrivial.}

The {\it central ratio} of a graph $G$ is the ratio of the cardinality of its center to its order; i.e., $|C(G)|/|V(G)|.$ In 1982, Buckley [1] proved the following result.

{\bf Theorem 1.} (Buckley [1]) {\it If $c$ is a rational number with $0<c\le 1,$ then there exists a graph whose central ratio is equal to $c.$}

In this paper, we obtain more detailed information by determining which cardinalities are possible for the center of a graph with given order and radius,
from which Theorem 1 follows immediately. There are unexpected phenomena in the results. We also prove a related uniqueness result.

\section{Main Results}

First recall the basic fact that if $n$ and $r$ are the order and radius of a graph respectively, then $r\le n/2.$ This can be seen by considering a spanning tree of the graph. It is well-known (e.g. [4]) that there exists a graph of radius $r$ and diameter $d$ if and only if
$r\le d\le 2r.$

We will need the following two lemmas due to other authors.

{\bf Lemma 2.} (Erd\H{o}s, Saks and S\'{o}s [2, Theorem 2.1])  {\it Every connected graph of radius $r$ has an induced path of order $2r-1.$}

A cycle $C$ in a graph $H$ is called a {\it geodesic cycle} if for any two vertices $x,y\in V(C),$ $d_{C}(x,y)=d_{H}(x,y).$ Thus, a geodesic cycle
is necessarily induced; i.e., it has no chords.

{\bf Lemma 3.} (Haviar, Hrn\v{c}iar and Monoszov\'{a} [3, Theorem 2.6]) {\it Let $G$ be a graph of order $n,$ radius $r$ and diameter $d.$
If $n\le 3r-2$ and $d\le 2r-2,$ then $G$ contains a geodesic cycle of length $2r$ or $2r+1.$}

 {\bf Notation 1.} For a positive integer $n,$\, $\langle n\rangle=\{1,2,\ldots,n\},$  the set of the first $n$ positive integers.

 $K_t,$ $C_t$ and $P_t$ will denote the complete graph of order $t,$ the cycle of order $t$ and the path of order $t$  respectively.
 Recall that the {\it join} of graphs $G$ and $H,$  denoted $G\vee H,$ is the graph obtained from the disjoint union $G+H$ by
 adding the edges $xy$ with $x\in V(G)$ and $y\in V(H).$ $\overline{G}$ denotes the complement of a graph $G.$
 A {\it leaf} in a graph is a vertex of degree $1.$

{\bf Notation 2.} A {\it lollipop} is the union of a cycle and a path having exactly one common vertex which is an end-vertex of the path.
We denote by $L(n,k)$ the lollipop of order $n$ whose cycle has length $k.$ A {\it broom} is a tree obtained by subdividing one edge of a star
an arbitrary number of times. We denote by $B(n,k)$ the broom of order $n$ and diameter $k.$ If $k\le n-2,$ the broom $B(n,k)$ has a unique vertex
of degree at least $3,$ which we call the {\it joint} of the broom. The {\it joint} of the path $B(n,n-1)=P_n$ is defined to be the neighbor of one of
the two end-vertices.

Two examples are depicted in Figure 1.
\vskip 3mm
\par
 \centerline{\includegraphics[width=4.8in]{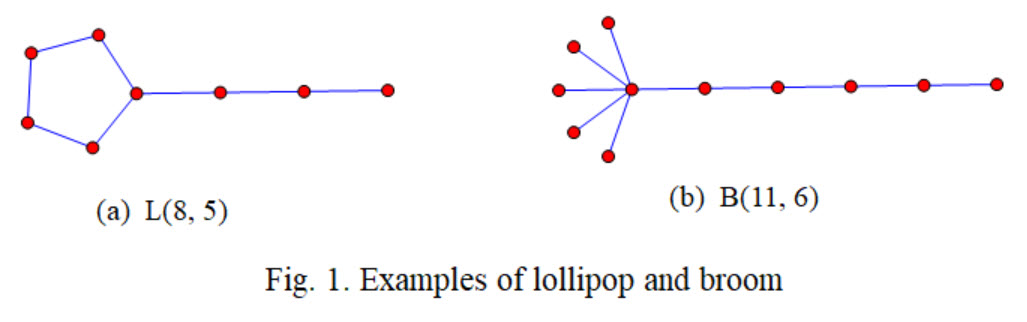}}
\par

 Now we are ready to state and prove the main results.

 {\bf Theorem 4.} {\it Given positive integers $n$ and $r$ with $r\le n/2,$\, let $\Omega(n,r)$ denote the set of those integers $k$
 such that there exists a graph of order $n$ and radius $r$ whose center has cardinality $k.$ If $n\ge 3$ then
$$
 \Omega(n,r)=\begin{cases} \langle n\rangle\setminus \{n-1\}\quad {\rm if}\,\,\,r\le (3n+2)/8;\\
                             \{1,2,...,n-2r+2\}\cup\{6r-2n+1,...,n-2,n\}\quad {\rm if}\,\,\,(3n+2)/8<r<n/2;\\
                             \{2,n\}\quad {\rm if}\,\,\,r=n/2.
                             \end{cases}
$$}

{\bf Proof.} We first show that every value $s$ in the expression of $\Omega(n,r)$ can be attained. Note that the condition $(3n+2)/8<r$
means $(n-2r+2)+1<6r-2n+1;$ i.e., there is an integer gap between $n-2r+2$ and $6r-2n+1.$

For this part, there is no need to distinguish the case $r\le (3n+2)/8$ and the case $(3n+2)/8<r<n/2.$ The constructions below give two graphs 
attaining some values of $s$ in the former case.

Suppose $r=1.$ For any $s\in \langle n\rangle\setminus \{n-1\},$ the graph $K_s\vee \overline{K_{n-s}}$ has order $n$ and radius $1,$
and its center has cardinality $s.$

Now suppose $2\le r<n/2.$ We distinguish six cases.

Case 1. $s=1.$ The broom $B(n,2r)$ has order $n$ and radius $r,$ and its center has cardinality $1.$

Case 2. $2\le s\le n-2r+2.$ Let $G_{1}(n,r,s)$ be the graph obtained from the broom $B(n-s+2,2r-1)$ by adding $s-2$ additional vertices
and joining each of them to the two central vertices of $B(n-s+2,2r-1).$ Then $G_{1}(n,r,s)$ has order $n$ and radius $r,$ and its center
has cardinality $s.$ The graph $G_{1}(14,4,5)$ is depicted in Figure 2.
\vskip 3mm
\par
 \centerline{\includegraphics[width=3 in]{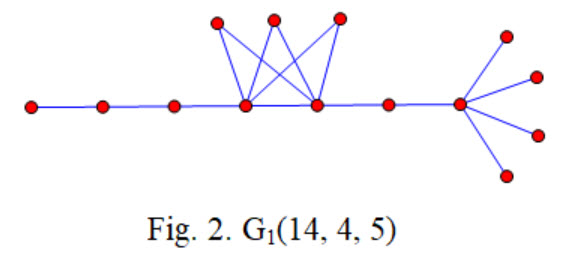}}
\par

Case 3. $6r-2n+1\le s\le 2r-1$ and $s$ is odd. Then $s=2r-2k+1$ for some $k$ with $1\le k\le n-2r.$
If $n=2r+k,$ define the graph $G_{2}(n,r,s)$ to be the lollipop $L(n,2r).$ If $n>2r+k,$ $G_{2}(n,r,s)$ is obtained
by identifying the joint of the broom $B(n-2r+1,k+1)$ with one vertex of the cycle $C_{2r}.$
Then $G_{2}(n,r,s)$ has order $n$ and radius $r,$ and its center has cardinality $s.$ The graph $G_{2}(15,4,3)$ is depicted in Figure 3.
\vskip 3mm
\par
 \centerline{\includegraphics[width=3.1 in]{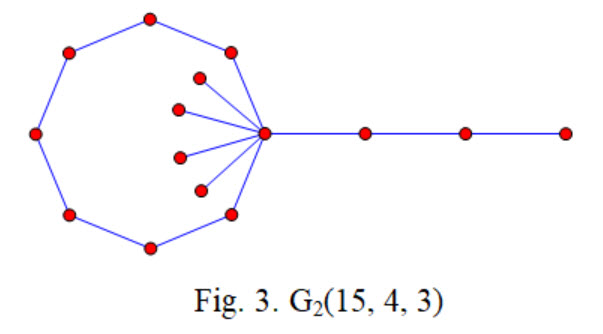}}
\par

Case 4. $6r-2n+1\le s\le 2r-1$ and $s$ is even. Then $s=2r-2k$ for some $k$ with $1\le k\le n-2r-1.$
If $n=2r+k+1,$ we identify the eccentric vertex of the unique leaf of the lollipop $L(2r+1,2r)$ with one end-vertex of
the path $P_{k+1}$ to obtain the graph $G_{3}(n,r,s).$ If $n>2r+k+1,$ we identify the eccentric vertex of the unique leaf of the lollipop $L(2r+1,2r)$
with the joint of the broom $B(n-2r,k+1)$ to obtain the graph $G_{3}(n,r,s).$ Then $G_{3}(n,r,s)$ has order $n$ and radius $r,$ and its center
has cardinality $s.$ The graph $G_{3}(15,4,2)$ is depicted in Figure 4.
\vskip 3mm
\par
 \centerline{\includegraphics[width=3.1 in]{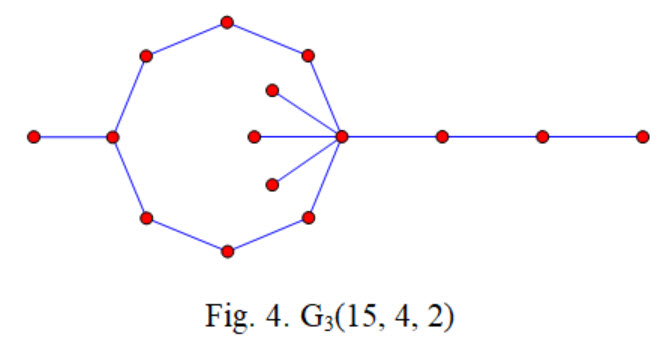}}
\par

Case 5. $2r\le s\le n-2.$ Then $s=2r+k-1$ for some $k$ with $1\le k\le n-2r-1.$ Define the graph $G_{4}(n,r,s)$ as follows.
Let $C:\,x_1,x_2,\ldots,x_{2r}$ be a cycle of length $2r.$ Add $n-2r-k$ new vertices to $C$ and join each of them to the vertex $x_1$ to
obtain a graph $D.$ Add $k$ new vertices to $D$ and join each of them to the two vertices $x_2$ and $x_{2r}$ to obtain the graph
$G_{4}(n,r,s),$ which has order $n$ and radius $r,$ and whose center has cardinality $s.$ The graph $G_{4}(15,4,11)$ is depicted in Figure 5.
\vskip 3mm
\par
 \centerline{\includegraphics[width=2 in]{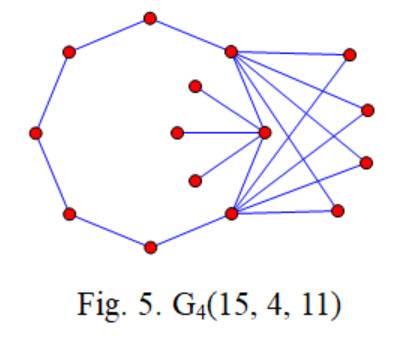}}
\par

Case 6. $s=n.$ We obtain a graph $G_{5}(n,r)$ by adding $n-2r$ new vertices to the cycle $C:\,x_1,x_2,\ldots,x_{2r}$ of length $2r$
and joining each of them to the two vertices $x_2$ and $x_{2r}.$ Then $G_{5}(n,r)$ is a self-centered graph of order $n$ and radius $r.$
The graph $G_{5}(12,4)$ is depicted in Figure 6.
\vskip 3mm
\par
 \centerline{\includegraphics[width=2 in]{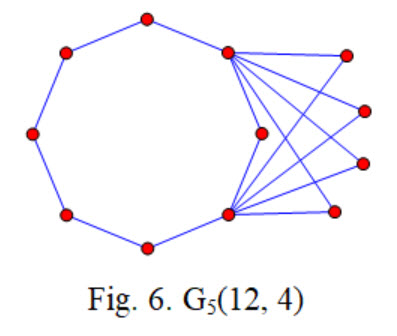}}
\par

Finally suppose $r=n/2.$ We have $n=2r.$ The center of the path $P_{2r}$ has cardinality $2$ and the center of the cycle $C_{2r}$
 has cardinality $n.$

Conversely, we prove that the integers appearing in the expression of $\Omega(n,r)$ are the only possible elements of the set $\Omega(n,r).$

First note that any nontrivial graph has at least two peripheral vertices. Indeed, suppose $v$ is a peripheral vertex and $u$
is an eccentric vertex of $v.$ Then both $v$ and $u$ are peripheral vertices. Thus, we always have $n-1\not\in \Omega(n,r).$

We then treat the easier case $r=n/2.$ Let $H$ be a graph of order $n$ and radius $r$ with $r=n/2.$ By Lemma 2, $H$
has an induced path $P$ of order $2r-1=n-1.$ Thus there is only one vertex $y$ outside $P.$ Since $P$ is an induced path of radius $r-1$
and $H$ has radius $r,$ the two end-vertices of $P$ are the only possible neighbors of $y.$ Hence $H$ is either a path or a cycle of the even order $n.$
Consequently the center of $H$ has cardinality $2$ or $n.$

Now we consider the case $(3n+2)/8<r<n/2.$ Let $G$ be a graph of order $n$ and radius $r$ whose center has cardinality $t.$ We will show that
either $t\le n-2r+2$ or $t\ge 6r-2n+1.$ Note that the condition $(3n+2)/8<r$ implies that $n\le 3r-2.$ Denote by $d$ the diameter of $G.$
We distinguish two cases.

Case 1. $d\ge 2r-1.$ Let $Q$ be a diametral path of $G;$ i.e., a path of length $d$ whose end-vertices are at distance $d.$
There are two possible values of $d.$ If $d=2r,$ then there are at least $2r$ non-central vertices of $G$ on $Q.$ Hence $t\le n-2r.$
If $d=2r-1,$ then there are at least $2r-2$ non-central vertices of $G$ on $Q.$ Hence $t\le n-2r+2.$ In both cases we have $t\le n-2r+2.$

Case 2. $d\le 2r-2.$  As remarked above, we also have $n\le 3r-2.$ By Lemma 3, $G$ has a cycle $C$ of length $2r$ or $2r+1.$
We distinguish two subcases.

Subcase 2.1. $C$ has length $2r.$ Since $2r<n,$ $G$ has vertices outside $C.$
Suppose $u_1v_1$ is a non-cut-edge of $G$ with $u_1\in V(C)$ and $v_1\not\in V(C).$ Recall that an edge is a
cut-edge if and only if it belongs to no cycle [5, p.23]. Thus there is a cycle containing $u_1v_1.$ In this cycle we choose an edge $e$ such that
$e\not= u_1v_1$ and $e\not\in E(C).$ Delete $e$ to obtain a new graph. If $u_1v_1$ is still a non-cut-edge of this new graph, delete an edge $f$ in a cycle
containing $u_1v_1$ such that $f\not= u_1v_1$ and $f\not\in E(C).$ After finitely many such edge deletion steps, we obtain a graph $G_1$ containing $C$
in which $u_1v_1$ is a cut-edge.

If $G_1$ contains a non-cut-edge $u_2v_2$ with exactly one endpoint in $C,$ by repeating the above edge deletion procedure,we obtain a graph
$G_2$ containing $C$ in which $u_2v_2$ is a cut-edge. Continuing this process, after finitely many steps we obtain a graph $R$ containing $C$
in which every edge with exactly one endpoint in $C$ is a cut-edge of $R.$ Suppose $u_iv_i,$ $i=1,2,\ldots, p$ are all such cut-edges of $R$
with $u_i\in V(C).$ The vertices $v_1,\ldots, v_p$ are pairwise distinct, but it is possible that $u_i=u_j$ for some pairs $i\neq j.$
Since deleting edges cannot decrease the radius of a graph, we have $r={\rm rad}(G)\le {\rm rad}(R).$ Denote by $F_i$ the component of $R-u_iv_i$ containing $v_i.$ Then $F_1,\ldots, F_p$ are pairwise vertex-disjoint.

Denote $a_i={\rm max}\{d_R(u_i,x)|x\in V(F_i)\},$ $i=1,\ldots, p.$ Since $C$ contains $2r$ vertices and $F_1,\ldots, F_p$ are pairwise vertex-disjoint, we have $a_1+a_2+\cdots+a_p\le n-2r.$ Since the order $n$ of $R$ satisfies $n\le 3r-2,$ we deduce that $a_i\le n-2r\le r-2$ for each $i.$ It follows that $C$ contains at most $2r-(2(r-a_i)+1)=2a_i-1$ vertices $y$ such that there exists a vertex $z\in V(F_i)$ with $d_R(y,z)>r.$ In $R,$ altogether $C$ contains at most
$$
\sum_{i=1}^p(2a_i-1)=2\left(\sum_{i=1}^p a_i\right)-p\le 2(n-2r)-1=2n-4r-1 \eqno (1)
$$
vertices with eccentricity larger than $r.$ Hence $C$ contains at least $2r-(2n-4r-1)=6r-2n+1$ vertices with eccentricity not exceeding $r$ in $R.$ Since
$6r-2n+1\ge 6r-2(3r-2)+1=5,$ we have ${\rm rad}(R)\le r.$ Combining this conclusion with the fact that ${\rm rad}(R)\ge r$ we obtain
${\rm rad}(R)=r;$ i.e., $R$ and $G$ have the same radius. Since $R$ is obtained from $G$ by deleting edges, every central vertex of
$R$ is also a central vertex of $G.$ We deduce that $R$ and hence $G$ contains at least $6r-2n+1$ central vertices.

Subcase 2.2. $C$ has length $2r+1.$ The proof is similar to that for Subcase 2.1, but there are some differences in details. Continue using
the above notations. $C$ contains at most $2\sum_{i=1}^p a_i$ vertices with eccentricity larger than $r.$ We have $2\sum_{i=1}^p a_i\le 2(n-2r-1).$
Hence, $C$ contains at least
$$
(2r+1)-2(n-2r-1)=6r-2n+3\ge 6r-2(3r-2)+3=7
$$
vertices of eccentricity $r.$ Thus, $G$ has at least $6r-2n+3$ central vertices.

Combining subcases 2.1 and 2.2 we conclude that $t\ge 6r-2n+1.$ This completes the proof. $\Box$

 {\bf Theorem 5.} {\it Let $n$ and $r$ be positive integers with $(3n+2)/8\le r<n/2.$ Then the lollipop $L(n,2r)$ is the unique graph
 of order $n$ and radius $r$ whose center has cardinality $6r-2n+1.$}

 {\bf Proof.} First note that the condition $(3n+2)/8\le r<n/2$ implies that $n\ge 7,$ $n\le 3r-2$ and $n-2r+2<6r-2n+1.$
 
 We observe that the lollipop $L(n,2r)$ has order $n$ and radius $r,$ and its center has cardinality $6r-2n+1.$

 Next suppose $G$ is a graph of order $n$ and radius $r$ whose center has cardinality $6r-2n+1.$
 Denote by $d$ the diameter of $G.$  We assert that $d\le 2r-2,$ since in the proof of Theorem 4 we have shown that 
 if $d\ge 2r-1,$ then the center of $G$ has cardinality at most $n-2r+2,$ which contradicts the fact that the center of $G$ has cardinality $6r-2n+1.$

 By Lemma 3, $G$ contains a geodesic cycle $C$ of length $2r$ or $2r+1.$ 
 In the proof of Theorem 4 we have shown that if $C$ has length $2r+1,$ then the center of $G$ has cardinality at least $6r-2n+3,$ which
 contradicts our assumption.  Hence $C$ has length $2r.$ Since $C$ is geodesic, it is an induced cycle. Continue using the notations
 in the proof of Theorem 4. Since  the center of $G$ has cardinality $6r-2n+1,$ equality must hold in (1), which implies that $p=1$
 and $\sum_{i=1}^p a_i=n-2r.$  Thus, in $R,$ the component $F_1$ is a path of order $n-2r,$ implying that $R$ is the lollipop $L(n,2r).$
 Note that the center of $L(n,2r)$ has cardinality $6r-2n+1.$  If in $G$ there is an edge other than $u_1v_1$ joining a vertex of $C$ and
 a vertex of $F_1,$ then the cardinality of the center of $G$ would be larger than $6r-2n+1,$ a contradiction.
 It follows that $G=R=L(n,2r)$ and the proof is complete. $\Box$

 For example, Theorem 5 asserts that the lollipop $L(14,12)$ is the unique graph of order $14$ and radius $6$ whose center has cardinality $9.$

 Finally we show that Buckley's Theorem 1 follows from Theorem 4.

 {\bf Proof of Theorem 1.} Any self-centered graph, say, a cycle, has central ratio $1.$ Now let $c=a/b$ with $a<b.$
 First suppose $a\neq b-1.$ Choose any positive integer $r$ with $r\le (3b+2)/8.$ By Theorem 4, there exists a graph $G$ of order $b$
 and radius $r$ whose center has cardinality $a.$ Then $G$ has central ratio $a/b.$ Next suppose $a=b-1.$ We have $2a\neq 2b-1.$
 Reasoning as above, there exists a graph $H$ of order $2b$ whose center has cardinality $2a.$ Clearly $H$ has central ratio
 $(2a)/(2b)=c.$ $\Box$

\vskip 5mm
{\bf Acknowledgement.} This research  was supported by the NSFC grants 11671148 and 11771148 and Science and Technology Commission of Shanghai Municipality (STCSM) grant 18dz2271000.

\end{document}